\documentclass[fleqn]{mat01}
\usepackage{times,mathtimy,amssymb,latexsym}

\newcommand{\bbC}{{\mathbb C}}
\newcommand{\bbZ}{{\mathbb Z}}

\newcommand{\bbP}{{\mathbb P}}
\newcommand{\bbF}{{\mathbb F}}

\newcommand{\bbQ}{{\mathbb Q}}
\newcommand{\Oh}{{\mathcal O}}

\newcommand{\calK}{{\mathcal K}}
\DeclareMathOperator{\Alb}{Alb}
\DeclareMathOperator{\Spec}{Spec}

\DeclareMathOperator{\HH}{\it H}

\DeclareMathOperator{\FF}{\it F}

\DeclareMathOperator{\KK}{\it K}
\DeclareMathOperator{\GG}{\it G}
\DeclareMathOperator{\CH}{\it CH}
\DeclareMathOperator{\Pic}{Pic}

\DeclareMathOperator{\PGL}{\it PGL}

\DeclareMathOperator{\deg0}{hom}
\DeclareMathOperator{\im}{Image}
\DeclareMathOperator{\coker}{coker}
\DeclareMathOperator{\SK}{\it SK}
\DeclareMathOperator{\reg}{reg}

\DeclareMathOperator{\supp}{supp}
\DeclareMathOperator{\ann}{Ann}
\newcommand{\onto}{\twoheadrightarrow}
\newcommand{\lonto}{\twoheadleftarrow}

\newcommand{\into}{\hookrightarrow}
\newcommand{\by}[1]{\xrightarrow{#1}}

\newcommand{\tensor}{\otimes}
\newcommand{\isom}{\cong}

\newenvironment{dgm}[1]{\arraycolsep=\doublerulesep\begin{array}{#1}
    }{\end{array}}

\begin{document}

\setcounter{page}{37}
\firstpage{37}

\font\xxxxx=msam10 at 10pt
\def\onto{\mbox{\xxxxx\,{\char'263}}\ \ \,}
\def\lonto{\mbox{\xxxxx\,{\char'264}}\ \ \,}

\newtheorem{lem}{Lemma}
\newtheorem{propo}{\rm PROPOSITION}
\newtheorem{theor}{\bf Theorem}
\newtheorem{conj}{Conjecture}
\newtheorem{cor}{\rm COROLLARY}
\newtheorem{sublem}{Sublemma}

\newtheorem{assum}{Assumption}
\newtheorem{defn}{\rm DEFINITION}

\newtheorem{remar}{Remark}
\newtheorem{quest}{Question}

\title{Zero cycles on certain surfaces in arbitrary characteristic}

\markboth{G~V~Ravindra}{Zero cycles on certain surfaces in
arbitrary characteristic}

\author{G~V~RAVINDRA}

\address{Department of Mathematics, Washington University, St.
Louis, MO 63130, USA\\
\noindent E-mail: ravindra@math.wustl.edu}

\volume{116}

\mon{February}

\parts{1}

\pubyear{2006}

\Date{MS received 6 October 2004; revised 5 November 2005}

\begin{abstract}
Let $k$ be a field of arbitrary characteristic. Let $S$ be a
singular surface defined over $k$ with multiple rational curve
singularities and suppose that the Chow group of zero cycles of
its normalisation $\tilde{S}$ is finite dimensional. We give
numerical conditions under which the Chow group of zero cycles of
$S$ is finite dimensional.
\end{abstract}


\maketitle

\section{Introduction}

This work arose out of an attempt to prove Bloch's conjecture for
the fake projective planes constructed by Mumford \cite{M} and
later by Ishida and Kato \cite{KI}. An earlier (again
unsuccessful) attempt to prove the same was by Barlow \cite{B}.

The strategy was to compute the Chow group of zero cycles of the
special fibre of a model of the fake planes (it turns out to be
isomorphic to $\bbZ$) and to use this to deduce that the generic
fibre, which is the fake plane, too has Chow group of zero cycles
isomorphic to $\bbZ$. Motivated by our example, we now pose (and
partly answer) the following general question:

\begin{quest}\label{aim}
{\rm Let $S$ be a projective surface defined over an algebraically
closed field such that $\tilde{S}$, the normalisation of $S$ is
smooth. Suppose now that $\CH_{0}(\tilde{S})$ is finite
dimensional (in a suitable sense). Then under what conditions is
$\CH_{0}(S)$ finite dimensional?}
\end{quest}

For instance, when $S$ is defined over the complex numbers it is
clear from Mumford's theorem \cite{M2} that the finite
dimensionality of the Chow group of zero cycles forces that there
are no non-trivial global two forms on $\tilde{S}$. So in this
case, one is interested in knowing what are the conditions under
which $S$ itself has a finite dimensional Chow group.

Our main theorem in this context is the following:
\begin{theor}[\!]\label{mainthm}
Let $S$ be a singular surface with multiple rational curve
singularities such that the normalisation $\tilde{S}$ is
non-singular. Further{\rm ,} assume that $\HH^{i}(S,\Oh_{S}) \isom
\HH^{i}(\tilde{S}, \Oh_{\tilde{S}})$ for $i=1,2$. Then if
$\CH_{0}(\tilde{S})$ is `finite dimensional'{\rm ,} then so is
$\CH_{0}(S)$.
\end{theor}

As an application we prove a consequence of Bloch's conjecture:

\begin{theor}[\!]
Let $M_0$ be the special fibre occurring in the construction of
Mumford's fake projective plane. Then $\CH_0(M_0)\isom \bbZ$.
\end{theor}

After this work was completed, we were informed that Gillet
\cite{HG} has proved similar results but for surfaces defined over
algebraically closed fields of characteristic zero. The
interesting aspect of our work is that we can obtain results for
fields of any characteristic. Further, we can show that the
cohomological conditions stated above imply a numerical condition
on the geometry of the surface which is very easy to check.

\section{Preliminaries}

\subsection{\it Definitions and generalities}

We work over a field $k$ of arbitrary characteristic.

For any variety $X$, let $\KK_{0}(X)$ (resp. $\GG_0(X)$) be the
Grothendieck group of vector bundles (resp. coherent sheaves). One
knows that when $X$ is smooth, $\GG_0(X)$ is isomorphic to
$\KK_0(X)$.

\begin{defn}$\left.\right.$\vspace{.5pc}

\noindent{\rm Let $S$ be an irreducible, reduced quasi-projective
surface over any field. Then
\begin{equation*}
\FF^{2}\KK_{0}(S):=\{\alpha \in \KK_{0}(S)~|~ {\rm
rank}(\alpha)=0= \det(\alpha)\}.
\end{equation*}}
\end{defn}

\begin{theor}[(Bloch-Levine)]
For $S$ as above{\rm ,} there is a natural isomorphism
\begin{equation*}
\CH_{0}(S) \isom \FF^{2}\KK_{0}(S).
\end{equation*}
\end{theor}

We now state Bloch's conjecture for surfaces, which is the main
interest of this work. First some generalities: For any smooth
projective variety $X$ of dimension $n$, there is a cycle class
map
\begin{equation*}
\CH_{n-p}(X) \by{cl_{X}} \HH^{2p}(X),
\end{equation*}
where the group on the right-hand side is a suitable Weil
cohomology group. Further, if we denote the kernel of this map by
$\CH_{n-p}(X)_{\deg0}$, then there is an Abel--Jacobi map for
$p=n$,
\begin{equation*}
\CH_0(X)_{\deg0} \by{AJ_{X}} \Alb_{X}.
\end{equation*}

\begin{defn}\label{fin-dim}$\left.\right.$\vspace{.5pc}

\noindent{\rm Let $X$ be a projective variety. $\CH_0(X)$ is said
to be {\em finite dimensional} (see \cite{M2}) if for some $m>0$,
the map
\begin{equation*}
\gamma_m\hbox{:}\ S^m(X_{\reg}) \to \CH_0(X)_{\deg0}
\end{equation*}
is surjective. Here $X_{\reg}$ refers to the smooth locus of the
variety $X$ and the map is the canonical one which after fixing a
base point $[x_0]$ sends an element $[x]$ to the class of the
difference $[x]-[x_0]$.}
\end{defn}

One can see that this is also equivalent to the statement that for
some integer $m'>0$, depending only on $X$, any element of
$\CH_0(X)_{\deg0}$ is represented by a 0-cycle
$\sum_{i=1}^r\delta_i$, where for each $i$, the cycle $\delta_i$
is a difference of two effective 0-cycles of degree $m'$ supported
in $X_{\reg}$.

For a smooth projective surface $S$ defined over an algebraically
closed field $k$, Roitman \cite{R} has shown that the Chow group
(of zero cycles) is finite dimensional precisely when the natural map
\begin{equation*}
\CH_0(S)_{\deg0} \to \Alb(S)
\end{equation*}
is an isomorphism.

\begin{conj}\hskip -.3pc{\rm ({\it Bloch}).}\ \ \
{\rm Let $S$ be a smooth projective surface defined over $\bbC$ with
$p_{g}=0$. Then $\CH_{0}(S)$ is finite dimensional.}
\end{conj}

The conjecture above has been proved when $S$ is not of general
type. Further, in characteristic $0$, the work of Gillet \cite{HG}
and later on the work of Krishna and Srinivas \cite{KS} sheds some
light on what controls the Chow group in the case when $S$ is
singular. Unfortunately, we do not know of any general conjecture
for singular surfaces defined over positive characteristic.

\subsection{\it The surfaces of interest}

Let $\tilde{S}$ be a smooth projective surface defined over $k$
and $Z'\subset \tilde{S}$ be a divisor. Suppose that there exists
a push out diagram
\begin{equation*}
\begin{dgm}{ccc}
Z' & \to & \tilde{S}\\
\downarrow & &\\
Z   &     &
\end{dgm},
\end{equation*}
where $Z' \to Z$ is a finite surjective morphism. Let
$\nu\hbox{:}~\tilde{S}\to S$ be the resulting morphism. Further
assume that the curves $Z'$ and $Z$ have the following properties:
\begin{enumerate}
\renewcommand\labelenumi{\arabic{enumi}.}
\leftskip -.1pc
\item $Z'$ is a union of smooth rational curves.

\item There are at most $3$ components passing through any point
in $Z$.

\item The components of $Z$ intersect transversally.

\item $Z'$ is a normal crossings divisor in $\tilde{S}$.
\end{enumerate}

Further, let $n_{1}$ be the number of components, $n_{2}$ the
number of points of two-fold intersections and $n_{3}$ the number
of points of three-fold intersections in $Z$. Similarly, let
$m_{1}$ (resp. $m_{2}$) be the number of components (resp. the
number of points of two-fold intersections) in $Z'$.

Here, by making a base change to a larger field, we can assume
without loss of generality that all points above are $k$-rational
points.

\begin{lem}\label{conducting}
In the situation above{\rm ,} $Z$ is the subscheme defined by the
conducting ideal.
\end{lem}

\begin{proof}
By definition we have to prove that $Z$ is the subscheme of $S$
defined by the $\Oh_{S}$-annihilator of the coherent sheaf
$\nu_{\ast}\Oh_{\tilde{S}}/\Oh_{S}$. Since $\mathcal{I}_Z=
\nu_{\ast}\mathcal{I}_{Z'}$, we have
\begin{equation*}
\nu_{\ast}\Oh_{\tilde{S}}/\Oh_{S}= \nu_{\ast}\Oh_{Z'}/\Oh_{Z}.
\end{equation*}
The annihilator of the latter is clearly $\mathcal{I}_{Z}$ and
thus we are done. (For any finitely generated module $M$, over a
Noetherian ring $A$, the zero set of the annihilator $\ann(M)$ is
equal to the support, $\supp({M})$.\hfill $\Box$
\end{proof}

\begin{prop}{\rm (Localisation sequence) \cite{KS,PW}.}
\label{localisationisom}$\left.\right.$\vspace{.5pc}

\noindent Let $(\tilde{S}, Z') \to (S, Z)$ be the normalisation
with $Z$ the subscheme defined by the conducting ideal and $Z'$
its inverse image. Then there exists a commutative diagram
\begin{equation}
\begin{dgm}{cccccccc}\label{localisation}
\to & \SK_1(Z) & \to & \FF^{2}\KK_{0}(S,Z) & \to & \FF^{2}\KK_0(S) & \to & 0\\
 &\downarrow  &  & \downarrow & &\downarrow  &  &\\
\to & \SK_1(Z') & \to & \FF^{2}\KK_{0}(\tilde{S},Z') & \to &
\FF^{2}\KK_0(\tilde{S}) &  \to & 0
\end{dgm},
\end{equation}
where for a curve $C$, $\SK_{1}(C)$ denotes the cohomology group
$\HH^{1}(C, \mathcal{K}_{2,C})$. Furthermore{\rm ,} in the diagram
above the middle vertical arrow is an isomorphism.
\end{prop}

\begin{cor}\label{thecorollary}$\left.\right.$\vspace{.5pc}

\noindent There is an exact sequence
\begin{equation*}
0 \to \coker(\SK_{1}(Z) \to \SK_{1}(Z')) \to \FF^{2}\KK_0(S) \to
\FF^{2}\KK_0(\tilde{S}) \to 0.
\end{equation*}
\end{cor}

\begin{proof}
Follows easily from diagram (\ref{localisation}) in
Proposition~\ref{localisationisom} above.\hfill $\Box$
\end{proof}

\begin{lem}\hskip -.3pc{\rm ({\it Mayer--Vietoris sequences}).}\ \ \
For the subschemes $Z$ and $Z'${\rm ,} there exist Mayer--Vietoris
sequences such that the following diagram commutes.
\begin{equation}
\hskip -4pc\begin{dgm}{ccccccccc} \calK_{2,Z} & \to &
\oplus_{i=1}^{n_{1}}\calK_{2,Z_i} & \by{\Phi} & \oplus_{1\leq i<j
\leq n_{1}}~ \calK_{2,Z_{ij}} & \to &
\oplus_{1\leq i<j<k \leq n_{1}}~ \calK_{2, Z_{ijk}} & \to & 0 \\
\downarrow{\nu_0^{\ast}}  &     & \downarrow{\nu_1^{\ast}}  &     &
\downarrow{\nu_2^{\ast}}  & &      & \\
\calK_{2,Z'} & \to & \oplus_{r=1}^{m_{1}}~ \calK_{2,Z_r'} & \by{\Phi'} &
 \oplus_{1\leq r<s \leq m_{1}}~ \calK_{2,Z_{rs}'} &
\to &  0  & & \\
\end{dgm}.
\end{equation}
\end{lem}

\begin{proof}
The Mayer--Vietoris sequence for $Z'$ (which is divisor of normal
crossings in $B$) can be obtained as follows: The map $\calK_{2,
Z'} \to \oplus_{i}\calK_{2, Z_i'}$ is the restriction map. The map
$\oplus_{r=1}^{m_{1}}~ \calK_{2,Z_r'} \to \oplus_{1\leq r<s \leq
m_{1}}~ \calK_{2,Z_{rs}'}$ is given by $(\alpha_{i}) \mapsto
(\phi_{i,j(i)}(\alpha_{i})-\phi_{i,j(i)}(\alpha_{j(i)}))$ where
for $i<j(i)$, $\phi_{i,j(i)}$ is the restriction map induced by
the inclusion of $Z_{ij(i)}'$ in $Z_{i}'$ (or $Z_{j(i)}'$). It is
enough to prove surjectivity of this map at the stalks. This
follows since $\KK_{2}(A)$ for a local ring $A$, is generated by
Steinberg symbols \cite{Mi} and hence the map on $\KK_{2}$ is
surjective for surjective maps of local rings.

The extra term in the Mayer--Vietoris sequence for the scheme $Z$
occurs because of the fact that the components do not necessarily
meet transversally. There are points (corresponding to some of the
two-fold intersections $Z_{ij}$) which occur more than once. Thus
one has to apply (Theorem~6.4 of \cite{Mi}) once more. Surjectivity now is
obvious.

Furthermore, the maps $\nu_1$ and $\nu_2$ are injective.

For the diagram to commute, one needs to define the map $\Phi'$
carefully: Notice that since the components of $Z$ are smooth,
this means that no two rational curves which get identified in $Z$
intersect in $Z'$. This implies that when we index the curves in
$Z'$, we need to index them according to the indexing in $Z$. For
instance, for two curves above $Z_r'$ and $Z_s'$, $r<s$ if and
only if their images $Z_i$ and $Z_j$ are such that $i<j$. Once
this is taken care of, it is easy to check that the diagram
commutes.\hfill $\Box$
\end{proof}

We rewrite the above diagram in a way which is suitable for
cohomology computations. Let
\begin{align*}
\overline\calK_{2,Z} &= \im\left(\calK_{2,Z} \to
\oplus_{i=1}^{n_{1}}~ \calK_{2,Z_i}\right),\\[.3pc]
\overline\calK_{2,Z'} &= \im\left(\calK_{2,Z'} \to
\oplus_{r=1}^{m_{1}}~ \calK_{2,Z_r'}\right).
\end{align*}

These give short exact sequences
\begin{align}
&\begin{dgm}{ccccccccc} 0 & \to & \overline\calK_{2,Z} & \to &
\oplus_{i=1}^{n_{1}}~ \calK_{2,Z_i}
 & \to &Q\phantom{00i} & \to & 0 \\
  &     &  \downarrow{\nu_0^{\ast}}  & & \downarrow{\nu_1^{\ast}}  &     &
\downarrow{\nu_2^{\ast}}  & &      \\
0 & \to & \overline\calK_{2,Z'} & \to &
\oplus_{r=1}^{m_{1}}~ \calK_{2,Z_r'} & \to &
\oplus_{1\leq r<s \leq m_{1}}~\calK_{2,Z_{rs}'} &
\to & 0 \\
\end{dgm},\label{mv1}\\[.4pc]
&0 \to Q  \to  \oplus_{1\leq i<j \leq n_{1}}~ \calK_{2,Z_{ij}} \to
\oplus_{1\leq i<j<k\leq n_{1}}~ \calK_{2, Z_{ijk}} \to
0.\label{mv2}
\end{align}

We record a couple of useful lemmas now.
\begin{lem}\label{tidbits}
For $Z'$ above{\rm ,} the following hold{\rm :}
\begin{enumerate}
\renewcommand\labelenumi{\rm \arabic{enumi}.}
\leftskip -.15pc
\item\label{nodal1} $\HH^{0}(Z_r',\calK_{2,Z'_r}) \isom \KK_2(k)$.

\item $\HH^{0}(Z_{rs}',\calK_{2,Z_{rs}'}) \isom \KK_2(k)^{\vert
Z_{rs}'\vert}$.

\item $\HH^{1}(Z',\overline\calK_{2,Z'}) \isom
\HH^{1}(Z',\calK_{2,Z'})$.

\item\label{nodal2} $\HH^{1}(Z_r',\calK_{2,Z'_r}) \isom k^{\ast}$.

\item $\HH^{1}(Z_{rs}',\calK_{2,Z_{rs}'})~=~0 $.\vspace{-.5pc}
\end{enumerate}
\end{lem}

\begin{proof}
Recall that $Z_r' \isom \bbP^1$ (this is assumption (1) in
\S~2.2).
\begin{enumerate}
\renewcommand\labelenumi{\rm \arabic{enumi}.}
\leftskip -.15pc
\item The result follows from the fact that $\KK_2(k[t])=\KK_2(k)$.

\item Obvious.

\item $\ker (\calK_{2,Z'} \onto \overline\calK_{2,Z'})$ is supported
on points.

\item Follows from the projective bundle formula and the BGQ spectral
sequence.

\item Obvious.\hfill $\Box$\vspace{-1.8pc}
\end{enumerate}
\end{proof}

Associated to diagram (\ref{mv1}), we have the following
commutative diagram of cohomology long exact sequences:
\begin{align}\label{cmv1}
\hskip -4pc
\begin{dgm}{c@{\,}c@{\,}c@{\,}c@{\,}c@{\,}c@{\,}c@{\,}c@{\,}c@{\,}c@{\,}c@{\,}c@{\,}c}
0 & \to & \HH^{0}(\overline\calK_{2,Z}) & \to &
\oplus\HH^{0}(\calK_{2,Z_i})
 & \by{\Phi} & \HH^{0}(Q) & \to & \HH^{1}(\calK_{2,Z}) & \to &
\oplus\HH^{1}(\calK_{2,Z_i}) & \to & 0 \\
  &     &  \downarrow{\epsilon_{1}}  & & \downarrow{\epsilon_{2}}  & &
\downarrow{\epsilon_{3}}  & & \downarrow{\epsilon_{4}}  & &
\downarrow{\epsilon_{5}}  & &      \\
0 & \to & \HH^{1}(\overline\calK_{2,Z'}) & \to &
\oplus\HH^{0}(\calK_{2,Z_r'}) & \by{\Phi'} &
\oplus\HH^{0}(\calK_{2,Z_{rs}'}) & \to &
\HH^{1}(\calK_{2,Z'}) & \to &
\oplus \HH^{1}(\calK_{2,Z_i'}) & \to & 0 \\
\end{dgm}
\end{align}

\begin{lem}
The cohomology long exact sequence associated to the sequence {\rm
(\ref{mv2})} yields $\HH^{0}(Q) \isom
\KK_{2}(k)^{\oplus{n_{2}-n_{3}}}$.
\end{lem}

\begin{proof}
Easily follows from the cohomology long exact sequence associated
to the short exact sequence (\ref{mv2}) and the fact that
$\HH^{1}(Q)=0$ as $Q$ is supported on points.\hfill $\Box$
\end{proof}

\begin{lem}\label{splitsurj}
Let $C$ be a projective curve defined over $k$ and $x \in C$ be a
$k$-rational point. Then the restriction map
$\HH^{0}(C,\mathcal{K}_{2,C}) \to \HH^{0}(\{x\},\mathcal{K}_{2,{x}})
\isom \KK_{2}(k)$ is a split surjection.
\end{lem}

\begin{lem}\label{fd}
Let $C$ be a smooth rational curve. Then
$\HH^{1}(C,\mathcal{K}_{2,C})$ is finite dimensional.
\end{lem}

\begin{proof}
This is immediate from the Gersten resolution for
$\mathcal{K}_{2,C}$. In fact one gets a surjection
$\Pic(C)\tensor{k^{\ast}} \onto \HH^{1}(C,\mathcal{K}_{2,C})$.
Since $\Pic(C)\isom \bbZ$, we are done.\hfill $\Box$
\end{proof}

\begin{lem}\label{answer}
Suppose $\CH_0(\tilde{S})$ is finite dimensional. Then $\CH_0(S)$
is finite dimensional if and only if the groups
$\im(\epsilon_{3})$ and $\im(\Phi')$ generate $\oplus_{1\leq r<s
\leq m_{1}}\HH^{0}(\calK_{2,Z_{rs}'})$ in diagram {\rm
(\ref{cmv1})}.
\end{lem}

\begin{proof}
$\coker(\epsilon_{5})$ is a quotient of $\oplus_{j}\HH^{1}
(\mathcal{K}_{2,Z_{j}'})$ and hence by Lemma~\ref{fd} (and some
diagram chasing) is finite dimensional. Thus the hypothesis in the
statement is equivalent to the condition that $\coker(\SK_{1}(Z)
\by{\epsilon_{4}} \SK_{1}(Z'))$ is finite dimensional.
Corollary~\ref{thecorollary} implies that under this hypothesis
\hbox{$\FF^{2}\!\!\KK_{0}(S)$} is finite dimensional if and only if
\hbox{$\FF^{2}\!\!\KK_{0}(\tilde{S})$} is so and thus we are done.\hfill
$\Box$
\end{proof}

\begin{theor}[\!]\label{mainresult}
The hypothesis in Lemma~{\rm \ref{answer}} holds if and only if
$m_{1}-m_{2} \geq n_{1}-n_{2}+n_{3}$.
\end{theor}

\begin{proof}
By Lemma~\ref{splitsurj}, diagram (\ref{cmv1}) reduces to the
following:\vspace{-.2pc}
\begin{equation}\label{inequality}
\begin{dgm}{ccccccccc}
  0 &  & 0         & &  & &  &  & \\
\downarrow & & \downarrow & &  & &  &  & \\
\KK_2(k)^{n_{1}} & \by{\Phi} & \KK_2(k)^{n_{2}-n_{3}} & \to & \SK_1(Z) & \to &
\oplus\HH^{1}(\calK_{2,Z_i}) & \to & 0 \\
\downarrow & & \downarrow{\epsilon_3} & & \downarrow{\epsilon_4}
 & & \downarrow{\epsilon_5} &  & \\
\KK_2(k)^{m_{1}} & \by{\Phi'} & \KK_2(k)^{m_{2}} & \to & \SK_1(Z') & \to &
\oplus \HH^{1}(\calK_{2,Z_i'}) & \to & 0 \\
\end{dgm}.
\end{equation}
Using the explicit description of the maps $\Phi$ and $\Phi'$ it
is not hard to check (we explicitly check these in the example in
the next section) that if $m_{1}-m_{2} \geq n_{1}-n_{2}+n_{3}$,
then the first vertical map is an isomorphism and thus
$\coker(\SK_{1}(Z) \to \SK_{1}(Z'))$ is isomorphic to
$\coker(\epsilon_{5})$ which as mentioned above is finite
dimensional. The converse is obvious.

\hfill $\Box$\vspace{-.5pc}
\end{proof}

The proof of theorem~\ref{mainthm} follows from the following.

\begin{prop}\label{ns}$\left.\right.$\vspace{.5pc}

\noindent Suppose that the natural maps $\HH^{1}(S,\Oh_{S}) \to
\HH^{1}(\tilde{S},\Oh_{\tilde{S}})$ and $\HH^{2}(S,\Oh_{S}) \to
\HH^{2}(\tilde{S},\Oh_{\tilde{S}})$ are isomorphisms. Then
$m_{1}-m_{2} \geq n_{1}-n_{2}+n_{3}$.
\end{prop}

\begin{proof}
Since $S$ is the pushout of the diagram\vspace{-.2pc}
\begin{equation*}
\begin{dgm}{ccc}
Z' & \to & \tilde{S}\\
\downarrow &    &\\
Z &   &\\[-.5pc]
\end{dgm}
\end{equation*}
we have seen earlier $\nu_{\ast}\mathcal{I}_{Z'} \isom
\mathcal{I}_{Z}$ and therefore\vspace{-.2pc}
\begin{equation*}
\frac{\nu_{\ast}\Oh_{\tilde{S}}}{\Oh_{S}} \isom
\frac{\nu_{\ast}\Oh_{Z'}}{\Oh_{Z}}.
\end{equation*}
Now consider the short exact sequence\vspace{-.2pc}
\begin{equation*}
0 \to \Oh_{S} \to \nu_{\ast}\Oh_{\tilde{S}} \to
\frac{\nu_{\ast}\Oh_{\tilde{S}}}{\Oh_{S}} \to 0.
\end{equation*}
The associated cohomology long exact sequence yields\vspace{-.2pc}
\begin{align*}
\hskip -4pc \cdots \HH^{1}(S, \Oh_{S}) \to \!\HH^{1}(\tilde{S},
\Oh_{\tilde{S}}) \to \!\HH^{1}
\left(\frac{\nu_{\ast}\Oh_{Z'}}{\Oh_{Z}}\right) \to \!\HH^{2}(S,
\Oh_{S}) \to \!\HH^{2}(S, \nu_{\ast}\Oh_{\tilde{S}}) \to 0.
\end{align*}

Here the right exactness follows from the fact that $\HH^{2}$ is
an exact functor in our situation ($S$ is a surface!!). The
hypothesis implies that\vspace{-.2pc}
\begin{equation*}
\HH^{1} \left(\frac{\nu_{\ast}\Oh_{\tilde{S}}}{\Oh_{S}}\right)
\isom \HH^{1} \left(\frac{\nu_{\ast}\Oh_{Z'}}{\Oh_{Z}}\right)= 0.
\end{equation*}
$\left.\right.$\vspace{-1pc}
\pagebreak

The sequence
\begin{equation*}
0 \to \Oh_{Z} \to \nu_{\ast}\Oh_{Z'} \to \frac{\nu_{\ast}
\Oh_{Z'}}{\Oh_{Z}} \to 0
\end{equation*}
yields a
surjection
\begin{equation*}
\HH^{1}(\Oh_{Z}) \to \HH^1(\nu_{\ast}\Oh_{Z'}).
\end{equation*}

To obtain the numerical condition, we consider the Mayer--Vietoris
sequences for the sheaves $\nu_{\ast}\Oh_{Z'}$ and $\Oh_{Z}$. We
imitate the computations done for comparing the sheaves
$\mathcal{K}_{2}$. For this, put
\begin{equation*}
Q':=\coker\{\Oh_{Z} \to \oplus_{i=1}^{n_1}\Oh_{Z_i}\}
\end{equation*}
to get
\begin{align*}
\begin{split}
&\begin{dgm}{ccccccccccc} 0 & \to & \Oh_{Z} & \to &
\oplus_{i=1}^{n_1}\Oh_{Z_i} & \to &
 Q'\phantom{00i} & \to & 0 \\
& & \downarrow{\nu_0^{\ast}} & & \downarrow{\nu_1^{\ast}} & &
\downarrow{\nu_2^{\ast}} & & & \\
0 & \to & \Oh_{Z'} & \to & \oplus_{1\leq i\leq m_1}\Oh_{Z_{i}'} & \to
& \oplus_{1\leq i<j \leq n_1}\Oh_{Z_{ij}'} & \to & 0 & &\\
\end{dgm}\\
&\begin{dgm}{ccccccccc} 0 & \to & Q' & \to & \oplus_{1\leq
i<j<n_1}\Oh_{Z_{ij}} & \to &
\oplus_{1\leq i<j<k<n_1}\Oh_{Z_{ijk}} & \to & 0 \\
\end{dgm}
\end{split}.
\end{align*}

Imitating the proof of theorem~\ref{mainresult} the result now
follows by noting that $\HH^{0}(Z,\Oh_{Z}) \isom
\HH^{0}(Z',\Oh_{Z'})\isom k$ and using the snake lemma to get a
short exact sequence
\begin{equation*}
0 \to \ker\{\HH^{1}(\Oh_{Z}) \to \HH^{1}(\nu_{\ast}\Oh_{Z'})\} \to
k^{m_1-n_1} \to k^{(m_2-n_2+n_3)} \to 0.
\end{equation*}

This implies that $(m_1-m_2)\geq (n_1-n_2+n_3)$.\hfill $\Box$
\end{proof}

\begin{remar}{\rm
In the case when $k=\bbC$, we note that if the conditions of
Proposition~\ref{ns} are satisfied then $h^2(\Oh_S)$ is zero if
and only if $h^2(\Oh_{\tilde{S}})$ is zero. Assuming Bloch's
conjecture, we get $\FF^{2}\KK_0(\tilde{S})$ is finite dimensional
and thus by Lemma~\ref{answer}, \hbox{$\FF^{2}\!\!\KK_0(S)$}, is also finite
dimensional.}
\end{remar}

\subsection{\it Representability questions}

This section is inspired by Theorem~B in \cite{HG}. By a theorem
of Roitman \cite{R}, it is known that if the Chow group of
$0$-cycles of a smooth projective surface over an algebraically
closed field $k$ is finite dimensional then it is representable.
In other words, there is an abelian variety, the Albanese of the
surface whose geometric points `represent' the degree $0$, $0$
cycles i.e., the (Albanese) map
\begin{equation*}
\CH_{0}(S)_{\deg0} \to \Alb(S)
\end{equation*}
is an isomorphism.

Theorem~B in \cite{Ku} says that in the case of complex projective
surfaces whose normalisation is smooth there is an isogeny
\begin{equation*}
\CH_{0}(S)_{\deg0} \to J^2(S),
\end{equation*}
where now $S$ is as above and $J^2(S)$ is the `naive' Albanese
i.e.,
\begin{equation*}
J^2(S) := \HH^{3}(S, \bbC)/\FF^{2}\!\!\HH^{3}(X,\bbC)+\HH^{3}(S,
\bbZ),
\end{equation*}
where $\FF^{2}\!\!\HH^{3}$ is the Hodge filtration defined by
Deligne.\pagebreak

Recently Esnault, Srinivas and Viehweg \cite{ESV} have proved the
following:

\begin{theor}[\!]\label{gap}
Let $X$ be a projective variety of dimension $n${\rm ,} defined
over an algebraically closed field $k$.
\begin{enumerate}
\renewcommand\labelenumi{\rm \arabic{enumi}.}
\leftskip -.15pc
\item There exists a smooth connected commutative algebraic group
$A^n(X)${\rm ,} together with a regular homomorphism
$\phi\hbox{\rm :}~\CH^n(X)_{\deg0}\to A^{n}(X)${\rm ,} such that
$\phi$ is universal among regular homomorphisms from the group
$\CH^n(X)_{\deg0}$ to smooth commutative algebraic groups.

\item Over a universal domain $k$ the Chow group is finite
dimensional precisely when $\phi$ is an isomorphism.

\item $A^n(X\times_kK)=A^n(X)\times_kK${\rm ,} for all algebraically
closed fields $K$ containing $k$.
\end{enumerate}
\end{theor}

%



The following now follows almost tautologically from the theorems
above.

\begin{theor}[\!]
Let $S$ be a projective surface such that the normalisation
$(\tilde{S}, Z') \to (S,Z)$ satisfies the conditions of
Proposition~{\rm \ref{ns}} above. Then if
$\CH_{0}(\tilde{S})_{\deg0}$ is finite dimensional {\rm (}i.e.{\rm
,} representable by $\Alb(\tilde{S})${\rm ),} then
$\CH_{0}(S)_{\deg0}$ is representable by $A^{2}(S)${\rm ,} the
generalised Albanese of Esnault{\rm ,} Srinivas and Viehweg.
\end{theor}

\begin{proof}
By Theorem~\ref{mainresult}, it follows that $\CH_{0}(S)_{\deg0}$
(which is \hbox{$\isom \FF^{2}\!\!\KK_{0}(S)$}) is finite dimensional. By
Theorem~\ref{gap} this implies that this group is isomorphic to
the generalised Albanese $A^2(S)$.\hfill $\Box$
\end{proof}

\section{An application: Bloch's conjecture for the special fibre}

\subsection{\it A brief description of the fake plane}

Let $R=\bbZ_{2}$ denote the ring of $2$-adic integers with
quotient field $K=\bbQ_{2}$ and finite residue field $k=\bbF_{2}$
of characteristic $2$. We describe here the fake plane constructed
by Mumford (henceforth denoted MFP) and refer the reader to
\cite{KI} for the other fake planes.

Mumford makes use of the method of $p$-adic uniformisation
introduced by Mustafin and Kurihara (see \cite{M} and the
references therein) to construct a formal scheme which is of
finite type over $\Spec(R)$. Take $\bbP^{2}_{R}$ and successively
blow up all the $k$-rational lines and $k$-rational points in the
special fibre. Let $U$ be the union of the generic fibre
$\bbP^2_K$ and a sufficiently small open neighbourhood of the
proper transform of $\bbP^2_k$ in the blown-up scheme above. For
each $\alpha \in \PGL(3,K)$, we denote by $U^{\alpha}$, the
$R$-scheme such that the generic fibre is equal to $\bbP^2_K$ and
that there exists an isomorphism $U \by{\isom} U^{\alpha}$ which
induces the natural action of $\alpha$ on the generic fibre. Then
the union of all $U^{\alpha}$ above is patched together to get a
regular scheme $X$. By construction, the action of $\PGL(3,K)$
extends to $X$. Let $\mathcal{X}$ be the completion of the
resulting scheme along the special fibre. For a certain discrete
torsion-free co-compact group $\Gamma \subset \PGL(3,K)$, he is
then able to construct a quotient subscheme $\mathcal{X}/\Gamma$.
Mumford ({\it op. cit.}) then shows that the canonical sheaf
$\omega_{\mathcal{X}}$ is ample and descends to
$\mathcal{X}/\Gamma$ and that the latter can be algebraized to a
projective variety denoted by $M$. For this choice of the subgroup
$\Gamma$, he shows that $\mathcal{X}/\Gamma$ has generic fibre, a
smooth surface of general type with the same Betti numbers as that
of the projective plane which we refer to as MFP. Further, the
special fibre is an irreducible rational surface over $k$ whose
normalisation is $\bbP^{2}_{k}$ blown-up at the $7$ rational
points.

Kato and Ishida \cite{KI} use results on discrete groups by
Cartwright {\it et~al} \cite{CMSZ} to construct two more fake
projective planes.

\subsection{\it The configurations of lines}

Below we give a description of the special fibre $M_{0}$ of the
fake planes. As noted above, the normalisation of the special
fibre, which we denote by $B$ is $\bbP^{2}_{k}$ blown up at its
$7$ rational points. The first and second tables give the
configuration of all rational curves in $B$ and the
identifications that give the configuration of lines in the
special fibre $M_{0}$ respectively and the last two tables give
the configuration of lines in the fake planes constructed by Kato
and Ishida. We will note here that while in the former one of the
rational double lines contains a node as a singularity, the latter
have no singular lines and the situation is as in the previous
situation.

Let $E(lmn)$ be the exceptional divisor over the point $[l:m:n]
\in \bbP^2_k$ and $C(lmn)$ be the proper transform of the line
given by $lx+my+nz=0$. We shall denote by $Z_{i}'$ and $Z_j$
respectively the rational curves in $B$ and $M_{0}$. Further, the
rational points in $M_{0}$ shall be denoted by $a$, $b$ etc. and
their preimages in $B$ shall also be denoted by the same\break
alphabet.

Tables~1 and 2 give the intersection data of the lines in the
configurations $Z'$ and $Z$.

\begin{table}[b]\vspace{.3pc} 
\processtable{Intersection data in the normalisation $Z'$}
{\begin{tabular}{@{}l@{\quad\ \ }c@{\quad\ \ }c@{\quad\ \ }c@{\quad\ \ }c@{\quad\ \ }c@{\quad\ \ }c@{\quad\ \ }c@{}}\hline\\[-.6pc]
     &$C(110)$&$C(100)$&$C(010)$&$C(001)$&$C(101)$&$C(011)$&$C(111)$\\
      &$=Z_8'$&$=Z_9'$&$=Z_{10}'$&$=Z_{11}'$&$=Z_{12}'$&$=Z_{13}'$&$=Z_{14}'$\\[.2pc]\hline\\[-.5pc]
$E(001)$ & $a$  & $b$  &  $a$    &         &         &         &        \\
$=Z_1'$  &      &      &         &         &         &         &         \\[.4pc]
$E(100)$&      &      &  $d$    &   $b$   &         &  $g$    &        \\
$=Z_2'$  &      &      &         &         &         &         &         \\[.4pc]
$E(110)$& $a$  &      &         &   $c$   &         &         &   $d$  \\
$=Z_3'$  &      &      &         &         &         &         &         \\[.4pc]
$E(111)$& $b$  &      &         &         &   $e$   &  $c$    &        \\
$=Z_4'$  &      &      &         &         &         &         &         \\[.4pc]
$E(011)$&      & $g$  &         &         &         &  $f$    &   $e$  \\
$=Z_5'$  &      &      &         &         &         &         &         \\[.4pc]
$E(101)$&      &      &  $c$    &         &   $g$   &         &   $f$  \\
$=Z_6'$  &      &      &         &         &         &         &         \\[.4pc]
$E(010)$&      & $d$  &         &   $e$   &   $f$   &         &        \\
$=Z_7'$  &      &      &         &         &         &         &\\[.2pc]\hline
\end{tabular}}{}
\end{table}

\begin{table}[b] 
\processtable{Intersection data in the special fibre $Z$}
{\begin{tabular}{@{}l@{\quad\ \ }c@{\quad\ \ }c@{\quad\ \ }c@{\quad\ \ }c@{\quad\ \ }c@{\quad\ \ }c@{\quad\ \ }c@{}}\hline\\[-.6pc]
     &$Z_1$ &$Z_2$ &  $Z_3$  & $Z_4$   &  $Z_5$  &  $Z_6$  &  $Z_7$\\[.1pc]\hline\\[-.6pc]
$Z_1$&      &      &  $e$    &  $g$    &  $f,~g$ & $e,~f$  &        \\[.2pc]
$Z_2$&      &      &  $c$    &   $d$   &  $c$    &  $d$    & $a,~a$ \\[.2pc]
$Z_3$& $e$  & $c$  &         &   $b$   &  $c$    &  $e$    &  $b$   \\[.2pc]
$Z_4$& $g$  & $d$  &  $b$    &         &  $g$    &  $d$    &  $b$   \\[.2pc]
$Z_5$&$f,~g$&$c$   &  $c$    &   $g$   &         &  $f$    &        \\[.2pc]
$Z_6$&$f,~e$& $d$  &  $e$    &   $d$   &  $f$    &         &        \\[.2pc]
$Z_7$&      &$a,~a$&  $b$    &   $b$   &         &         &  $a$\\[.2pc]\hline
\end{tabular}}{}\vspace{-.8pc}
\end{table}

Note that every exceptional curve meets a proper transform of a
rational line in $Z$ in each of its rational points ($3$ in number)
and similarly a proper transform meets an exceptional line in each of
its rational points (again $3$). Further no two exceptional or proper
transform intersect.

In the above one can easily check that $Z_i'$ is identified with
$Z_{i+7}'$ to obtain the configuration $Z$.

The following are the intersection tables (tables~3, 4) of the
lines in the special fibres of fake planes constructed by
Kato--Ishida.

\begin{table}[b]\vspace{.7pc} 
\processtable{} {\begin{tabular}{@{}l@{\quad\ \ }c@{\quad\ \ }c@{\quad\ \ }c@{\quad\ \ }c@{\quad\ \ }c@{\quad\ \ }c@{\quad\ \ }c@{}}\hline\\[-.6pc]
      &$Z_8'$&$Z_9'$&$Z_{10}'$&$Z_{11}'$&$Z_{12}'$&$Z_{13}'$&$Z_{14}'$\\[.15pc]\hline\\[-.6pc]
$Z_1'$&      &      &  $f$    &         &  $e$    &         &  $g$   \\[.4pc]
$Z_2'$&      &      &         &   $d$   &         &  $c$    &  $b$   \\[.4pc]
$Z_3'$& $g$  &      &         &   $a$   &  $f$    &         &        \\[.4pc]
$Z_4'$&      & $c$  &         &         &  $a$    &         &  $d$   \\[.4pc]
$Z_5'$& $f$  &      &  $a$    &         &         &  $e$    &        \\[.4pc]
$Z_6'$& $e$  & $b$  &     &       $c$   &         &         &        \\[.4pc]
$Z_7'$&      & $d$  &  $g$    &         &         &  $b$    &        \\[.3pc]\hline
\end{tabular}}{}
\end{table}

\begin{table} 
\processtable{} {\begin{tabular}{@{}l@{\quad\ \ }c@{\quad\ \ }c@{\quad\ \ }c@{\quad\ \ }c@{\quad\ \ }c@{\quad\ \ }c@{\quad\ \ }c@{}}\hline\\[-.6pc]
      &$Z_8'$&$Z_9'$&$Z_{10}'$&$Z_{11}'$&$Z_{12}'$&$Z_{13}'$&$Z_{14}'$\\[.15pc]\hline\\[-.6pc]
$Z_1'$&      &      &  $e$    &         &   $g$   &         &   $f$  \\[.4pc]
$Z_2'$&      & $b$  &  $b$    &         &         &  $c$    &        \\[.4pc]
$Z_3'$&      & $b$  &         &   $a$   &   $e$   &         &        \\[.4pc]
$Z_4'$& $g$  &      &         &         &   $a$   &  $d$    &        \\[.4pc]
$Z_5'$& $e$  &      &  $a$    &   $g$   &         &         &        \\[.4pc]
$Z_6'$&      &      &         &   $d$   &         &  $d$    &   $c$  \\[.4pc]
$Z_7'$& $f$  & $c$  &         &         &         &
&$f$\\[.3pc]\hline
\end{tabular}}{}\vspace{-.7pc}
\end{table}

\subsection{\it Precise statement of the conjecture}

Let $F$ be a universal domain containing $k$, $W$ the Witt vectors
of $F$ and let $M$ denote the base change of the original $M$
(over $\bbZ_2$) to $W$. Let $M_0$ and $M_{\eta}$ be the special
fibre and the generic fibre of $M$ defined over $F$ (which is now
an infinite extension of $k$) and $\Omega$ (which denotes the
quotient field of $W$) respectively. Let $B \to M_0$ be the
normalisation. One knows from above that $B$ is isomorphic to the
base change to $F$ of $\bbP^2_k$ blown up at all $k$-rational
points.

By Riemann--Roch one has $\GG_0(M_0)\tensor\bbQ \isom \oplus
\CH_i(M_0)\tensor\bbQ$. Since $B \to M_0$ is a finite map, one has
a surjection $\GG_0(B)\tensor\bbQ \to \GG_0(M_0)\tensor\bbQ$.
Since the first group is finite rank $\bbQ$-vector space (for
instance, by the Grothendieck--Riemann--Roch theorem) so is the
second.

We note the following important consequence of Bloch's conjecture for
fake projective planes.

\begin{lem}
If Bloch's conjecture is true for the fake planes{\rm ,} then the
group \hbox{$\FF^{2}\!\!\KK_0(M_0)\tensor\bbQ$} is a finite dimensional
$\bbQ$-vector space.
\end{lem}

\begin{proof}
Consider the following diagram:
\begin{equation*}
\begin{dgm}{ccccccc}
\GG_0(M_0)\tensor\bbQ & \to & \GG_0(M)\tensor\bbQ
 & \to & \GG_0(M_{\eta})\tensor\bbQ & \to & 0 \\
           &     & \downarrow &     &                 &     &    \\
           &     &  \KK_0(M_0)\tensor\bbQ &    & & & \\
           &     &  \downarrow &    & & & \\
           &     &       0      &    & & &
\end{dgm}.
\end{equation*}
Here the horizontal sequence is the localisation theorem for the
$K$-theory of (coherent sheaves on) $M$ and the vertical map is
the restriction of $K$-groups. (Since $M$ is regular, in this case
we have $\GG_0(M) \isom \KK_0(M)$.) The surjectivity of the
vertical map follows from the fact that any line bundle on $M_0$
lifts since $p_g=0=q$ for $M_0$ and since any element of
$\KK_0(M_0)$ is equivalent to a sum of line bundles. If Bloch's
conjecture holds for $M_{\eta}$ then this implies that
$\GG_0(M_{\eta})\tensor\bbQ \left(\isom \oplus
\CH_i(M_{\eta})\tensor\bbQ \right) \isom \bbQ^{3}$. Further from
the previous paragraph, we know that $\GG_0(M_0)\tensor\bbQ$ is of
finite rank.\hfill $\Box$
\end{proof}

\subsection{\it The Chow group of the special fibre}

In this section, we compute the Chow group of the special fibre
$M_0$ of the fake projective planes. Ishida \cite{I} and
Kato--Ishida \cite{KI} have explicitly described the
identifications of the seven pairs of lines $C(a,b,c)$ and
$E(a,b,c)$ of the normalisation $B$. We use these identifications
in an essential way in our computation of the Chow group
\hbox{$\FF^{2}\!\!\KK_{0}(M_0)$}. We compute the above group in the case when
$M$ is Mumford's fake projective plane and leave the fake
projective planes of Kato and Ishida as an exercise to the reader.

Let $Z=\cup_{i=1}^{7}Z_{i}$ be the union of the $F$-rational lines
$Z_i$ in the special fibre $M_0$ and
$Z_{i}^{'}=\cup_{i=1}^{14}Z_{i}^{'}$ be the union of the $7$ pairs
of $F$-rational lines in the normalisation $B$ of $M_0$ (namely
the $7$ exceptional divisors and the $7$ proper transforms of the
lines in $\bbP^{2}_F$).

From the previous section we have for the map $\nu_0\hbox{:}~Z'
\to Z$ (the restriction of the normalisation $\nu\hbox{:}~B \to
M_0$), a diagram of $\calK$-sheaves:
\begin{equation}\label{me}
\begin{dgm}{c@{\,}c@{\,}ccccccc}
\calK_{2,Z} & \to & \oplus_{i=1}^{7}\calK_{2,Z_i} & \to &
\oplus_{1\leq i<j \leq 7}~ \calK_{2,Z_{ij}} & \to &
\oplus_{1\leq i<j<k \leq 7}~ \calK_{2, Z_{ijk}} & \to & 0 \\
\downarrow{\nu_0^{\ast}}  &     & \downarrow{\nu_1^{\ast}}  &     &
\downarrow{\nu_2^{\ast}}  & &      & \\
\calK_{2,Z'} & \to & \oplus_{r=1}^{14}~ \calK_{2,Z_r'} & \to &
 \oplus_{1\leq r<s \leq 14}~ \calK_{2,Z_{rs}'} &
\to & 0  & &
\end{dgm}
\end{equation}

Using the intersection table, it is easy to check that the above
diagram is indeed commutative.  From the earlier section, we have
short exact sequences
\begin{align}
&\begin{dgm}{ccccccccc} 0 & \to & \overline\calK_{2,Z} & \to &
\oplus_{i=1}^{7}~ \calK_{2,Z_i} & \to
& Q\phantom{00i} & \to & 0 \\
  &     &  \downarrow{\nu_0^{\ast}}  & & \downarrow{\nu_1^{\ast}}  &     &
\downarrow{\nu_2^{\ast}}  & &      \\
0 & \to & \overline\calK_{2,Z'} & \to & \oplus_{r=1}^{14}~
\calK_{2,Z_r'} & \to & \oplus_{1\leq r<s \leq
14}~\calK_{2,Z_{rs}'} & \to & 0
\end{dgm},\label{me1}\\[.4pc]
&0 \to Q  \to  \oplus_{1\leq i<j \leq 7}~ \calK_{2,Z_{ij}} \to
\oplus_{1\leq i<j<k\leq 7}~ \calK_{2, Z_{ijk}} \to   0.\label{me2}
\end{align}

Similar statements as in Lemma~\ref{tidbits} hold for the variety $Z$
but in this case one needs the following lemma.

\begin{lem}\label{nodaltidbits}
Let $C$ be a nodal rational curve. Then there are surjections
$\HH^{i}(\mathcal{K}_{2,C}) \to \HH^{i}(\mathcal{K}_{2,\bbP^1})$ for
$i=0,1$.
\end{lem}
\vspace{-1pc}
\begin{proof}
The desingularisation $\bbP^1 \to C$ induces a sequence
\begin{equation*}
0 \to K \to \mathcal{K}_{2,C} \to \mathcal{K}_{2,\bbP^1} \to L \to
0.
\end{equation*}
Here the extreme two terms are supported on the singular point.
The cohomology sequence associated to this yields a sequence
\begin{align*}
\hskip -4pc 0 \to \HH^0(K) \to \HH^{0}(\mathcal{K}_{2,C}) \by{f}
\HH^0(\mathcal{K}_{2,\bbP^1}) \to\! \HH^0(L) \to\!
\HH^{1}(\mathcal{K}_{2,C}) \to\! \HH^1(\mathcal{K}_{2,\bbP^1}) \to 0.
\end{align*}
Since $C$ has a point, this implies that the map $f$ splits. Thus
we see that the maps $\HH^{i}(\mathcal{K}_{2,C}) \to
\HH^{i}(\mathcal{K}_{2,\bbP^1})$ are surjective for $i=0,1$.\hfill $\Box$
\end{proof}

Consider the exact sequence (\ref{me2}). Taking the long exact
sequence of cohomology one~gets
\begin{align*}
0 &\to \HH^{0}(Q) \to \oplus_{1\leq i<j \leq 7}~
\HH^{0}(\calK_{2,Z_{ij}}) \to \oplus_{1\leq i<j<k \leq 7}~
\HH^{0}(\calK_{2,Z_{ijk}})\\[.3pc]
&\to \HH^{1}(Q) \to \oplus_{1\leq i<j\leq 7}~ \HH^{1}(\calK_{2,Z_{ij}})
\to 0.
\end{align*}

Since $Q$ is supported on points, $\HH^{1}(Q)$ is zero. Using the
previous lemma, one gets
\begin{equation*}
0 \to \HH^{0}(Q) \to \KK_2(F)^{\oplus 20} \to \KK_2(F)^{\oplus 6}
\to 0.
\end{equation*}
This implies that $\HH^{0}(Q) \isom \KK_2(F)^{\oplus 14}$.
\vspace{-.1pc}
\begin{lem}
$\HH^{0}(Z,\overline\calK_{2,Z}) \isom
\HH^{0}(Z',\overline\calK_{2,Z'}) \isom \KK_2(F)$.
\end{lem}
\vspace{-1pc}
\begin{proof}
That $\HH^{0}(Z',\overline\calK_{2,Z'}) \isom \KK_2(F)$ follows from
an inspection of the cohomology long exact sequence of the
Mayer--Vietoris sequence of $Z'$.

Note that the map
\begin{equation*}
\oplus_{1\leq r \leq 14}~\HH^{0}(\calK_{2,Z_{r}'}) \to
\oplus_{1\leq r<s \leq 14}~\HH^{0}(\calK_{2,Z_{rs}'})
\end{equation*}
is abstractly the map $ \KK_2(F)^{\oplus 14} \to \KK_2(F)^{\oplus
21}$ which takes the $r$th component of an element
$(\underline\alpha) \in \KK_2(F)^{\oplus 14}$, $\alpha_{r}\mapsto
(\beta_{ij})$ where $\beta_{rs}=\alpha_{r}-\alpha_{s}$ for $s$
such that $Z_{s}'$ intersects $Z_{r}'$ and $\beta_{ij}=0$
otherwise. Following table $1$, it is easy to see that the kernel
of this map is isomorphic to $\KK_2(F)$.

To see that a similar statement holds true for $Z$, one observes
that the map $\HH^{0}(Q) \to \oplus_{1\leq r<s \leq
14}\HH^{0}(\calK_{2,Z_{rs}})$ is injective. Using the
identifications that follow from tables $1$ and $2$ and the
explicit description of maps in the Mayer--Vietoris sequences, one
can then check that the kernel of the composite map
\begin{equation*}
\KK_2(F)^{\oplus 7} \to \KK_2(F)^{\oplus 14} \to \KK_2(F)^{\oplus
21}
\end{equation*}
is isomorphic to $\KK_2(F)$.

This implies that $\HH^{0}(Z,\overline\calK_{2,Z}) \isom \KK_2(F)$.\hfill $\Box$
\end{proof}

Taking the long exact sequence of cohomology for the diagram
(\ref{me1}) one gets
\begin{align}\label{cohom1}
\begin{dgm}{c@{\ }c@{\ }c@{\ }c@{\ }c@{\ }c@{\ }c@{\ }c@{\ }c@{\ }c@{\ }c@{\ }c@{\ }c}
0 & \to & \KK_2(F) & \to &
\KK_2(F)^{\oplus 7}  & \to
& \KK_2(F)^{\oplus 14} & \to &
\SK_1(Z) & \to & (F^{\ast})^{\oplus 7} & \to &  0 \\
  &     & \downarrow &   &
\downarrow           &
& \downarrow           &     &
\downarrow &   & \downarrow            & & \\
0 & \to & \KK_2(F) & \to &
\KK_2(F)^{\oplus 14}  & \to
& \KK_2(F)^{\oplus 21} & \to &
\SK_1(Z') & \to & (F^{\ast})^{\oplus 14} & \to & 0\\
\end{dgm}
\end{align}

\begin{remar}{\rm
Here one of the components say $Z_0$ is a nodal curve. We replace the
group $\HH^0(\mathcal{K}_{2,Z_0})$ by its image in $\KK_2(F)$. By
Lemma~\ref{nodaltidbits}, this is $\KK_2(F)$.}
\end{remar}

Using the earlier computations we get a commutative diagram
\begin{equation}\label{generation}
\begin{dgm}{ccccccccc}
0 & \to & \KK_2(F)^{\oplus 8} & \to & \SK_1(Z) & \to &
(F^{\ast})^{\oplus 7} & \to & 0 \\
& &\downarrow{\isom} & & \downarrow & & \downarrow &  & \\
0 & \to & \KK_2(F)^{\oplus 8} & \to & \SK_1(Z') & \to &
(F^{\ast})^{\oplus 14} & \to & 0
\end{dgm}.
\end{equation}

\begin{lem}\label{cokerofsk1}
\!$\SK_{\!1}(Z')$ is generated by $\SK_{\!1}(Z)$ and $\SK_{\!1}(B)$. In
particular{\rm ,} $\coker(\SK_{\!1}(Z)\!\!\to$ $\SK_{\!1}(Z'))$ is finite
dimensional.
\end{lem}

\begin{proof}
To understand the image of the map $\SK_1(B) \to
\SK_1(Z')$ we consider the following diagram:
\begin{equation*}
\begin{dgm}{ccccc}
\SK_1(B) & \onto & \HH^{1}(B,\calK_{2,B}) & \lonto & \Pic(B)\tensor\KK_1(F) \\
         &       &  \downarrow            &        &    \\
         &        & \SK_1(Z') & &
\end{dgm}.
\end{equation*}
Thus the $\im(\SK_1(B) \to \SK_1(Z'))= \im(\Pic(B)\tensor\KK_1(F)\to
\SK_1(Z'))$. From diagram (5), it is clear that all we need to check
is that $\oplus_{r}\HH^{1}(\calK_{2,Z_r'})$ is generated by the images
of $\oplus_{i}\HH^{1}(\calK_{2,Z_i})$ and $\SK_1(B)$. Abstractly this
is just saying that a free $F^{\ast}$-module with generators
$\{z_r\}_{r=1}^{14}$ is generated by the elements
$\{z_r+z_{r+7}\}_{r=1}^7$ and $\{z_r\}_{r=1}^7$ which is indeed
obvious.\hfill $\Box$
\end{proof}

\begin{lem}\label{middleterm}
In the localisation sequence {\rm (}diagram {\rm (\ref{localisation}))} for the map
$B \to M_{0}${\rm ,} the vertical map \hbox{$\FF^{2}\!\!\KK_{0}(M_0,Z) \to
\FF^{2}\!\!\KK_{0}(B,Z')$} is an isomorphism.
\end{lem}

\begin{proof}
Since $\nu\hbox{:}~B \to M_0$ is the normalisation, the above statement
follows from Proposition~\ref{localisationisom}, if one checks that $Z
\subset M_0$ is the subscheme defined by the conducting subscheme for
the normalisation morphism $\nu$. Since $M_{0}$ is obtained as a
pushout of the diagram
\begin{equation*}
\begin{dgm}{ccc}
Z' & \into  & B \\
\downarrow & & \\
Z  & &
\end{dgm},
\end{equation*}
this follows from the same argument as in Lemma~\ref{conducting}.\hfill $\Box$
\end{proof}

Now we come to the main result of this section.

\begin{theor}[\!]
The Chow group $\FF^{2}\!\!\KK_0(M_0)$ is isomorphic to $\bbZ$.
\end{theor}

\begin{proof}
The result follows by a trivial diagram chase in the localisation
sequence below:
\begin{equation}\label{f.d.}
\begin{dgm}{ccccccccc}
 \SK_1(Z) & \to & \FF^{2}\!\!\KK_{0}(M_0,Z)
& \to & \FF^{2}\!\!\KK_0(M_0) & \to & 0 \\
\downarrow  &  & \downarrow
&  &\downarrow  &  &  \\
 \SK_1(Z') & \to & \FF^{2}\!\!\KK_{0}(B,Z')
& \to & \FF^{2}\!\!\KK_0(B) &  \to & 0
\end{dgm}.
\end{equation}


By Lemma~\ref{cokerofsk1}, $\coker(\SK_1(Z) \to \SK_1(Z'))$ is
finite dimensional. Moreover the middle vertical arrow is an
isomorphism by Lemma~\ref{middleterm}. A~simple diagram chase will
then show that \hbox{$\FF^{2}\!\!\KK_0(M_0)$} is isomorphic to
\hbox{$\FF^{2}\!\!\KK_0(B)$} up to a finitely generated free
$F^{\ast}$-module (= $\coker(\SK_1(Z) \to \SK_1(Z'))$). Since $B$
is a smooth rational surface, \hbox{$\FF^{2}\!\!\KK_0(B)\isom\bbZ$}. Thus
\hbox{$\FF^2\!\!\KK_0(M_0)$} being a finitely generated abelian group, it is
finite dimensional. Now by Theorem~\ref{gap}, the Albanese map is
an isomorphism. Since this generalised Albanese is connected this
implies that $\FF^2\!\!\KK_0(M_0)\isom\bbZ$.\hfill $\Box$
\end{proof}

\begin{remar}{\rm
We direct the attention of the reader interested in further
examples to the paper of Gillet \cite{HG} mentioned in the
Introduction.}
\end{remar}

\section*{Acknowledgements}

It gives us great pleasure to acknowledge Najmuddin Fakhruddin for
suggesting the problem and for all the ideas that he generously
shared with us. In particular, the crucial idea that one can use
the localisation sequence to compute the Chow group of zero cycles
is due to him. We are also grateful to him for carefully reading
the manuscript and his suggestions. We would like to acknowledge
the support and encouragement we have received from V~Srinivas
over the past few years. We also thank him for clarifying certain
issues regarding the work \cite{ESV}. We thank A~Krishna for
bringing to our notice the article by Gillet~\cite{HG}.

\end{document}